\documentclass[a4paper,11pt]{article}
\usepackage{amsfonts,amssymb,amsmath}

\usepackage[english]{babel}
 \makeatletter
 \@addtoreset{equation}{section}
 \makeatother

 \makeatletter
 \@addtoreset{enunciato}{section}
 \makeatother

 \newcounter{enunciato}[section]

 \newtheorem{ittheorem}{Theorem}
 \newtheorem{itlemma}{Lemma}
 \newtheorem{itproposition}{Proposition}
 \newtheorem{itcorollary}{Corollary}
 \newtheorem{itdefinition}{Definition}
 \newtheorem{itremark}{Remark}
 \newtheorem{itclaim}{Claim}
 \newtheorem{itfact}{Fact}
 \newtheorem{itconjecture}{Conjecture}

 \newenvironment{theorem}{\addtocounter{enunciato}{1}
 \begin{ittheorem}}{\end{ittheorem}}

 \newenvironment{lemma}{\addtocounter{enunciato}{1}
 \begin{itlemma}}{\end{itlemma}}

 \newenvironment{proposition}{\addtocounter{enunciato}{1}
 \begin{itproposition}}{\end{itproposition}}

 \newenvironment{corollary}{\addtocounter{enunciato}{1}
 \begin{itcorollary}}{\end{itcorollary}}

 \newenvironment{definition}{\addtocounter{enunciato}{1}
 \begin{itdefinition}}{\end{itdefinition}}

 \newenvironment{remark}{\addtocounter{enunciato}{1}
 \begin{itremark}}{\end{itremark}}

 \newenvironment{claim}{\addtocounter{enunciato}{1}
 \begin{itclaim}}{\end{itclaim}}

 \newenvironment{fact}{\addtocounter{enunciato}{1}
 \begin{itfact}}{\end{itfact}}

 \newenvironment{conjecture}{\addtocounter{enunciato}{1}
 \begin{itconjecture}}{\end{itconjecture}}

 \newcommand{\be}[1]{\begin{equation}\label{#1}}
 \newcommand{\ee}{\end{equation}}

 \newcommand{\bl}[1]{\begin{lemma}\label{#1}}
 \newcommand{\el}{\end{lemma}}

 \newcommand{\br}[1]{\begin{remark}\label{#1}}
 \newcommand{\er}{\end{remark}}

 \newcommand{\bt}[1]{\begin{theorem}\label{#1}}
 \newcommand{\et}{\end{theorem}}

 \newcommand{\bd}[1]{\begin{definition}\label{#1}}
 \newcommand{\ed}{\end{definition}}

 \newcommand{\bcl}[1]{\begin{claim}\label{#1}}
 \newcommand{\ecl}{\end{claim}}

 \newcommand{\bfact}[1]{\begin{fact}\label{#1}}
 \newcommand{\efact}{\end{fact}}

 \newcommand{\bp}[1]{\begin{proposition}\label{#1}}
 \newcommand{\ep}{\end{proposition}}

 \newcommand{\bc}[1]{\begin{corollary}\label{#1}}
 \newcommand{\ec}{\end{corollary}}

 \newcommand{\bcj}[1]{\begin{conjecture}\label{#1}}
 \newcommand{\ecj}{\end{conjecture}}

 \newcommand{\bpr}{\begin{proof}}
 \newcommand{\epr}{\end{proof}}

 \newcommand{\bprlem}[1]{\begin{proofof}{\it\ref{#1}}.\,\,}
 \newcommand{\eprlem}{\end{proofof}}

 \newcommand{\bi}{\begin{itemize}}
 \newcommand{\ei}{\end{itemize}}

 \newcommand{\ben}{\begin{enumerate}}
 \newcommand{\een}{\end{enumerate}}

 \newcommand{\reff}[1]{(\ref{#1})}

 \newenvironment{proof}{\noindent {\em Proof}.\,\,}{\hspace*{\fill}$\halmos$\medskip}
 \newenvironment{proofof}{\noindent {\em Proof of Lemma\,\,}}{\hspace*{\fill}$\halmos$\medskip}
 \newcommand{\halmos}{\rule{1ex}{1.4ex}}
 \def \qed {{\hspace*{\fill}$\halmos$\medskip}}

 \newcommand{\one}{{\mathchoice {1\mskip-4mu\mathrm l}
         {1\mskip-4mu\mathrm l}
         {1\mskip-4.5mu\mathrm l}
         {1\mskip-5mu\mathrm l}}}

\def \qed {{\hspace*{\fill}$\halmos$\medskip}}

\def \L {{\mathbb L}}

\def \R {{\mathbb R}}
\def \Z {{\mathbb Z}}

\def \ra {\rightarrow}
\def \ba {\begin{array}}
\def \ea {\end{array}}

\def \Id {{\rm Id}}

\def \var {{\rm var}}
\def \osc {{\rm osc}}

\def \c {{\rm c}}

\def \cA {{\mathcal A}}

\def \cE {{\mathcal E}}
\def \cF {{\mathcal F}}
\def \cG {{\mathcal G}}
\def \cH {{\mathcal H}}

\def \cS {{\mathcal S}}
\def \cSb {{\mathcal S}_b}

\def\one{\rlap{\mbox{\small\rm 1}}\kern.15em 1}

\def\embf#1{\emph{\bf #1}}

\begin{document}
\title{Uniqueness and non-uniqueness of chains on half lines}

\author{\renewcommand{\thefootnote}{\arabic{footnote}}
R.\ Fern\'andez
\footnotemark[1]
\\
\renewcommand{\thefootnote}{\arabic{footnote}}
G.\ Maillard
\footnotemark[2]
}

\footnotetext[1]{
Laboratoire de Math\'ematiques Rapha\"el Salem,
UMR 6085 CNRS-Universit\'e de Rouen,
Avenue de l'Universit\'e, BP 12,
F-76801 Saint \'Etienne du Rouvray, France,
{\sl roberto.fernandez@univ-rouen.fr}
}
\footnotetext[2]{
Institut de Math\'ematiques, \'Ecole Polytechnique F\'ed\'erale
de Lausanne, CH-1015 Lausanne, Switzerland,
{\sl gregory.maillard@epfl.ch}
}

\maketitle

\begin{abstract}
We establish a one-to-one correspondence between one-sided and two-sided regular
systems of conditional probabilities on the half-line that preserves the associated
chains and Gibbs measures. As an application, we determine uniqueness and non-uniqueness
regimes in one-sided versions of ferromagnetic
Ising models with long range interactions. Our study shows that the interplay between
chain and Gibbsian theories yields more information than that contained
within the known theory of each separate framework.  In particular: (i) A Gibbsian
construction due to Dyson yields a new family of chains with phase transitions;
(ii) these transitions show that a square
summability uniqueness condition of chains is false in the general
non-shift-invariant setting, and (iii) an uniqueness criterion for chains shows that
a Gibbsian conjecture due to Kac and Thompson is false in this half-line setting.

{\it Acknowledgment.} The authors wish to thank Sacha Friedli and
Charles-Edouard Pfister for many enlightening discussions.
\end{abstract}

\newpage


\section{Introduction and preliminaries}
\label{S1}


\subsection{Introduction}
\label{S1.1}

Non-Markovian processes bring in the novel future of phase transitions: Several measures can
share the same transition probabilities if these have a sufficiently strong dependence on
faraway future~\cite{brakal93,berhofsid05,hul06}.   Unlike the Markovian case, this coexistence
is not due to the partition of the space into non-communicating components ---the transition
probabilities are all strictly positive in the published examples--- but rather to the
persistence of the influence of past history into the infinite future.  Such transitions
parallel statistical mechanical first-order phase transitions, where different boundary
conditions lead, in the thermodynamic limit, to different consistent measures.  This suggests,
as advocated in~\cite{fermai05}, to take Gibbsian theory as a model for the study of multiple-chain
phase diagrams.

More directly, one may wonder whether chains can simply be treated as one-dimensional Gibbs measures.
If so, the usual theory of discrete-time processes ---geared towards the description of phenomena
characteristic of Markov processes--- could be supplemented by an appropriately transcribed Gibbs
theory ---tailored to the description of complicated phase diagrams. Our attempt in~\cite{fermai04}
was unable to reach the multiple-phase case.

From a complementary point of view, it is natural to search for conditions granting that phase
transitions do not occur, that is, granting that a given family of transition probabilities admit
only one consistent measure.  There exist, at present, a number of such uniqueness
criteria~\cite{har55,ber87,lal86,lal00,comferfer02,sten03,fermai05,hul06,johobe03,johobepol05}
involving different non-nullness and continuity hypotheses and yielding information on different
properties of the invariant measure (mixing properties, Markovian approximation schemes, regeneration
and perfect simulation procedures).

In all these studies, there is an ignored aspect that deserves, in our opinion,
more careful consideration:  the role of the shift-invariance \emph{of the transition probabilities}.
All the phase transition examples involve shift-invariant transitions (and measures) and this is also
an ingrained feature in most of the proofs of existing uniqueness criteria.  The only exceptions are
the ``regeneration'' criterion of~\cite{comferfer02} and the ``one-side bounded uniformity'' proven
in~\cite{fermai05}. (For notational simplicity, a translation invariant setting was adopted
in~\cite{comferfer02}, but it is clear that the proof ---showing that every finite window can be
reconstructed from a finite past--- does not require shift invariance.)  It is legitimate to inquire
whether shift-invariance is an unavoidable requirement for the remaining criteria, or only an artifact
of the proof.

In this paper we illustrate some of the differences and similarities that arise when shift invariance is
lost.  We consider chains defined on the half line or, equivalent, on a ``time'' axis that is a countable
set with a total order and a maximal element.  We show, first, that in this setting
we can successfully complete the program initiated in~\cite{fermai04} and establish
a full correspondence between chains and Gibbs measures (Theorem \ref{th1} below).
We exploit the interplay
between Gibbsian and chain points of view to reveal a number of interesting facts:

\begin{itemize}
\item[(i)]
We borrow results by Dyson~\cite{dys69a} to show that the chains defined by
long range Ising models with couplings
decaying as power laws $|i-j|^{-p}$ ($p>1$ to ensure summability of the interactions) exhibit phase
transitions for $p<2$ for large values of the coupling parameters (``low temperatures'', see Theorem \ref{ising-res}).
The transcription of Dyson's approach amounts to a novel way to prove phase transitions in the context of chains,
namely by constructing a measure that is not mixing.  This implies that there should be at least two different
extremal consistent measures~\cite{fermai05}.
\item[(ii)]
The chains with $3/2<p<2$ do satisfy the square summability condition of Johansson and \"Oberg's uniqueness
criterion~\cite{johobe03} (see also~\cite{johobepol05}) and yet exhibit phase transitions (Remark~\ref{johobe-rk}
below).  This shows that such criterion ---which has been proven to be optimal in an appropriate
sense~\cite{berhofsid05}--- is false in general non-shift-invariant settings.
\item[(iii)]
Using the regeneration (chain) criterium of~\cite{comferfer02} we prove that, at least in the half line, a
conjecture by Kac and Thompson (mentioned in~\cite{dys69a}) is false (see Remark~\ref{ruel-rem} below).
\item[(iv)]
On the other hand, Gibbsian uniqueness criteria can be used to show that these models have a unique invariant
state at high temperatures.  The only chain criterion that is temperature-sensitive is one-sided
Dobrushin~\cite{fermai05,hul05} which, however, is not directly applicable to the Ising chains considered here.
\item[(v)] Present work implies a one-sided version of the Kozlov theorem~\cite{koz74} (Corollary \ref{koz-chains}
below): transition probabilities that are continuous and non-null are always defined by one-sided interactions,
albeit in an indirect manner that passes through an auxiliary specification.
\end{itemize}


\subsection{Notation and preliminary definitions}
\label{S1.2}

We consider a measurable space $(\cA,\cE)$ where $\cA$ is a finite alphabet and $\cE$ is the
discrete $\sigma$-algebra. We denote $(\Omega,\cF)$ the associated product measurable space
with $\Omega=\cA^{\L}$, where $\L$ is a countable set with total order. In this paper we study
the case in which $\L$ has a maximal element, that is, $\L=\Z^-$, but we shall also refer to
the usual unbounded case, where $\L=\Z$.  For each $\Lambda\subset\L$ we denote
$\Omega_\Lambda = \cA^{\Lambda}$ and $\sigma_\Lambda$ for the restriction of a configuration
$\sigma\in\Omega$ to $\Omega_\Lambda$, namely the family
$(\sigma_i)_{i\in\Lambda}\in \cA^{\Lambda}$.  Also, $\cF_\Lambda$ will denote the
sub-$\sigma$-algebra of $\mathcal{F}$ generated by cylinders based on $\Lambda$
($\cF_\Lambda$-measurable functions are insensitive to configuration values outside $\Lambda$).
When $\Lambda$ is an interval, $\Lambda=[k,n]$ with $k,n\in\L$ such that $k\le n$, we use the
notation: $l_\Lambda=k$, $m_\Lambda=n$, $\Lambda_-=\{i\in\L:i<k\}$,
$\omega_{k}^{n}=\omega_{[k,n]}=\omega_{k},\ldots ,\omega_{n}$, $\Omega_k^n=\Omega_{[k,n]}$ and
$\cF_k^n=\cF_{[k,n]}$. For semi-intervals we denote also $\cF_{\le n}:=\cF_{(-\infty,n]}$, etc.
The concatenation notation $\omega_\Lambda\,\sigma_\Delta$, where $\Lambda\cap\Delta=\emptyset$,
indicates the configuration on $\Lambda\cup\Delta$ coinciding with $\omega_i$ for $i\in\Lambda$
and with $\sigma_i$ for $i\in\Delta$. We denote $\cS$ the set of finite subsets of $\L$ and $\cSb$
the set of finite intervals of $\L$. To lighten up formulas involving probability kernels, we will
freely use $\nu(h)$ instead of $E_\nu(h)$ for $\nu$ a measure on $\Omega$ and $h$ a
$\mathcal{F}$-measurable function. Also $\nu(\sigma_\Lambda)$ will mean
$\nu(\{\omega\in\Omega:\omega_\Lambda=\sigma_\Lambda\})$ for
$\Lambda\subset\L$ and $\sigma_\Lambda\in\Omega_\Lambda$.

For any sub-$\sigma$-algebra $\cH$ of $\cF$, we recall that a \emph{measure kernel} on $\cH\times\Omega$
is a map $\pi(\,\cdot\mid\cdot\,)\colon$ $\cH\times\Omega\to\mathbb{R}$ such that $\pi(\,\cdot\mid \omega)$
is a measure on $(\Omega,\cH)$ for each $\omega\in\Omega$ while $\pi(A\mid\cdot\,)$ is $\cF$-measurable
for each event $A\in\cH$.  If each $\pi(\,\cdot\mid \omega)$ is a probability measure the kernel is called
a \emph{probability kernel}.
For kernels
$\pi$ and $\widetilde\pi$, non-negative measurable functions $h$, measures $\nu$ on $(\Omega,\cH)$
and cylinders $C_{\omega_\Lambda}=\{\sigma\in\Omega: \sigma_\Lambda=\omega_\Lambda\}$, we shall denote:

\begin{itemize}
\item $\pi(\omega_\Lambda\mid\,\cdot\,)$ for $\pi(C_{\omega_\Lambda}\mid\,\cdot\,)$
\item $\pi(h)$ for the measurable function $\int_\Omega
  h(\eta)\,\pi(d\eta\mid\,\cdot\,)$.
\item $\pi\widetilde\pi$ for the composed kernel defined by
  $(\pi\widetilde\pi)(h) = \pi\big(\widetilde\pi(h)\big)$.
\item $\nu\pi$ for measure defined by $(\nu\pi)(h) =\nu\big(\pi(h)\big)$.
\end{itemize}

In the unbounded case, $\L=\Z$, the (right) \emph{shift operator}
---$\tau:\Omega\to\Omega$, $(\tau \omega)_i=\omega_{i-1}$--- is an isomorphism that naturally induces
shift operations for measurable functions and kernels: $(\tau f)(\omega) = f(\tau^{-1}\omega)$,
$(\tau\pi)(h\mid\omega) = \pi(\tau^{-1}f\mid\tau^{-1}\omega)$.


\subsection{Chains and Gibbs measures}
\label{S1.2.2}\label{S1.2.1}

We start by briefly reviewing in parallel the well known notions of chains and Gibbs measures
in the spirit of~\cite{fermai05}.  The main difference is that, in chains, kernels apply only
to functions measurable with respect to the present and the past.

\bd{lis1}
A \embf{left singleton-specification (LSS)} (or system of transition probabilities) $f$
on $(\Omega,\cF)$ is a family of probability kernels $\left\{f_i\right\}_{i\in\Z}$ with
$f_i\colon\cF_{\leq i}\times\Omega\to [0,1]$ such that for all $i\in\Z$,
\begin{itemize}
\item[\rm{(a)}] for each $A\in\cF_{\leq i}$, $f_i(A\mid\cdot\,)$ is
$\cF_{\leq i-1}$-measurable;
\item[\rm{(b)}] for each $B\in\cF_{\leq i-1}$ and $\omega\in\Omega$,
$f_i(B \mid \omega)=\one_B(\omega)$.
\end{itemize}
The LSS $f$ is:
\begin{itemize}
\item[\rm{(i)}]
\embf{Continuous} if the functions
$f_i\left(\omega_i \mid \cdot\,\right)$ are continuous
for each $i \in \L$ and $\omega_i\in\Omega_i$;
\item[\rm{(ii)}]
\embf{Non-null} if the functions
$f_i\left(\omega_i \mid \cdot\,\right)$ are (strictly) positives for each
$i \in \Z$ and $\omega_i\in\Omega_i$;
\item[\rm{(iii)}]
\embf{Regular} if it is continuous and non-null;
\item[\rm{(iv)}]
\embf{Shift-invariant} if $\L=\Z$ and $\tau f_i=f_{i+1}$, $i\in\L$.
\end{itemize}
\ed
\bd{chain}
A probability measure $\mu$ on $(\Omega,\; \mathcal{F} )$ is said to be \embf{consistent}
with a LSS $f$ if for each $i \in\Z$,
\be{}
\mu f_i = \mu \quad\text{over }
\cF_{\leq i}\;.
\ee
The family of these measures will be denoted by $\mathcal{G}(f)$ and each $\mu\in\mathcal{G}(f)$
is called a {($f$-) \embf{chain}}. A measure $\mu$ is a \embf{regular chain} if there exists a
regular LSS $f$ such that $\mu\in\cG(f)$.
\ed

The singletons $f_i$ of a LSS define, through compositions, interval-kernels
\be{lss-k}
f_{[m,n]} \;=\; f_m\,f_{m+1}\,\cdots\, f_m
\ee
for $m\le n\in\L$.  We observe that
\be{lis-mu}
\mu\in\mathcal{G}(f) \;\Longleftrightarrow\; \mu f_{[m,n]} = \mu
\quad\text{over }\cF_{\leq n}\quad \forall m\le n\in\L\;.
\ee
The family $\{f_{[m,n]} : m\le n\in\L\}$ ---called a LIS (\emph{Left Interval-Specification})
in~\cite{fermai04,fermai05}--- is the chain analogue of the notion of specification.

\br{g-function}
The particular case when $f$ and $\mu\in\cG(f)$ are shift-invariant, reduces
to the study of $g$-functions and $g$-measures, respectively~\cite{kea72}.
Chains for general, non-shift-invariant singletons have also been called
$G$-measures~\cite{brodoo91,brodoo98}.
\er
\bd{specif}
A \embf{specification} $\gamma$ on $\left(\Omega,\cF\right)$ is a family of probability
kernels $\left\{\gamma_\Lambda\right\}_{\Lambda\in\cS}$ with
$\gamma_\Lambda:\cF\times\Omega\rightarrow [0,1]$ such that for all $\Lambda$ in $\cS$,
\begin{itemize}
\item[\rm{(a)}] for each $A\in\cF$, $\gamma_\Lambda(A\mid\cdot)$ is $\cF_{\Lambda^{c}}$-measurable;
\item[\rm{(b)}] for each $B\in\cF_{\Lambda^{c}}$ and $\omega\in\Omega$,
$\gamma_\Lambda(B\mid\omega)=\one_B(\omega)$;
\item[\rm{(c)}] for each $\Delta\in\cS:\Delta\supset\Lambda$,
$\gamma_\Delta \gamma_\Lambda=\gamma_\Delta$.
\end{itemize}
A specification $\gamma$ is:
\begin{itemize}
\item[\rm{(i)}]
\embf{Continuous} if the functions
$\gamma_\Lambda\left(\omega_\Lambda \mid \cdot\,\right)$ are continuous
for each $\Lambda \in \cS$ and $\omega_\Lambda\in\Omega_\Lambda$;
\item[\rm{(ii)}]
\embf{Non-null} if the functions
$\gamma_\Lambda\left(\omega_\Lambda \mid \cdot\,\right)$ are (strictly) positives
for each $\Lambda \in \cS$ and $\omega_\Lambda\in\Omega_\Lambda$;
\item[\rm{(iii)}]
\embf{Gibbsian} if it is continuous and non-null;
\item[\rm{(iv)}]
\embf{Shift-invariant} if $\L=\Z$ and $\tau \gamma_\Lambda=\gamma_{\Lambda+1}$, $\Lambda\in\cS$.
\end{itemize}
\ed
\bd{consist}
A probability measure $\mu$ on $(\Omega, \mathcal{F})$ is said to be \embf{consistent}
with a specification $\gamma$ if for each $\Lambda \in \cS$,
\be{gibbs9}
\mu \,\gamma_{\Lambda}\;=\;\mu.
\ee
The family of these measures will be denoted by $\cG(\gamma)$.
A measure $\mu$ is a \embf{Gibbs measure} if there exists a Gibbsian specification $\gamma$
such that $\mu\in\cG(\gamma)$.
\ed

A celebrated theorem due to Kozlov~\cite{koz74} (see also~\cite{sul73} for an alternative
version in a different interaction space) shows that a specification $\gamma$ is Gibbsian if,
and only if, there exists a \embf{potential}, i.e.\ a family of functions
$\phi=(\phi_A)_{A\in\cS}$ with each $\phi_A:\Omega\ra\R$ being $\cF_A$-measurable, that is
absolutely and uniformly summable in the sense
\be{abssum}
\sum_{A\ni i}\|\phi_A\|_\infty<\infty\,\,\forall\,i\in\Z;
\ee
such that $\gamma=\gamma^\phi$ where for all $\Lambda\in\cS$ and $\omega,\sigma\in\Omega$
\be{gibbsspe}
\gamma_\Lambda^\phi(\sigma_\Lambda\mid\omega_{\Lambda^c})
= \frac{\exp\big[-H_{\Lambda,\omega}^\phi(\sigma)\big]}{Z_{\Lambda,\omega}^\phi}\;,
\ee
 with
\be{gibbsspe*}
H_{\Lambda,\omega}^\phi(\sigma)
=\sum_{{A\in\cS}\atop{A\cap\Lambda\neq\emptyset}}\phi_A(\sigma_\Lambda\,\omega_{\Lambda^\c})
\quad\text{and}\quad
Z_{\Lambda,\omega}^\phi
=\sum_{\sigma_\Lambda\in\Omega_\Lambda}\exp\Big[-H_{\Lambda,\omega}^\phi(\sigma)\Big].
\ee
This is the original statistical mechanical prescription due to Boltzmann and Gibbs.


\section{Main Results}
\label{S2}

Let $\Omega=\cA^{\L}$ with $\L=\Z^-:=\{\cdots,-2,-1,0\}$, be the configuration space on the half-line ending at $0$.


\subsection{Correspondence between regular LSS and Gibbsian specifications on half-spaces}
\label{S2.1}

Let
\be{1to1-sets}
\begin{aligned}
\Theta &= \Big\{\text{regular LSS on } \Omega\Big\}\\
\Pi &= \Big\{\text{Gibbsian specifications on } \Omega\Big\}\\
\end{aligned}
\ee
and introduce the following maps
\be{1to1-maps}
b\colon\Theta\ra\Pi,\,\, f\mapsto\gamma^f
\quad\text{and}\quad
c\colon\Pi\ra\Theta,\,\, \gamma\mapsto f^\gamma
\ee
defined by
\be{2.2.1}
\gamma^f_{[l,0]}\;=\; f_{[l,0]}\;;
\ee
\be{1to1-maps*}
\gamma_\Lambda^f(\sigma_\Lambda\mid\omega_{\Lambda^\c})
=\frac{f_{[l_\Lambda,0]}\big(\sigma_\Lambda\,\omega_{[l_\Lambda,0]\cap\Lambda^\c}\mid \omega_{\Lambda_-}\big)}
{f_{[l_\Lambda,0]}\big(\omega_{[l_\Lambda,0]\cap\Lambda^\c}\mid \omega_{\Lambda_-}\big)}
\quad\text{with}
\qquad\omega\in\Omega\;,\,\Lambda\in\cS\,\;,
\ee
with $\Lambda$ strictly contained in $[l_\Lambda,0]$, and
\be{1to1-maps**}
f_i^\gamma\big(\omega_i\mid\omega_{-\infty}^{i-1}\big)
=\gamma_{[i,0]}\big(\omega_i\mid\omega_{-\infty}^{i-1}\big)
\qquad\forall\,i\in\Z,\,\omega\in\Omega\;.
\ee
We observe that, due to the consistency of $\gamma$,
\be{2.4.1}
f^\gamma_{[l,0]} \;=\; \gamma_{[l,0]}\;.
\ee
\bt{th1}
The maps $b$ and $c$ establish a one-to-one correspondence between $\Theta$ and $\Pi$
that preserves consistency. More precisely,
\begin{itemize}
\item[{\rm 1)}]
\begin{itemize}
\item[{\rm (a)}] $f^\gamma\in\Theta$;
\item[{\rm (b)}] $\gamma^f\in\Pi$;
\end{itemize}
\item[{\rm 2)}]
\begin{itemize}
\item[{\rm (a)}] $b\circ c=\Id_\Pi$;
\item[{\rm (b)}] $c\circ b=\Id_\Theta$;
\item[{\rm (c)}] $\cG\big(\gamma^f\big)=\cG(f)$, and
\item[{\rm (d)}] $\cG\big(f^\gamma\big)=\cG(\gamma)$.
\end{itemize}
\end{itemize}
\et

\bc{koz-chains}
For every regular LSS $f$ there exists an absolutely and uniformly summable potential
$\phi$ such that, for each $i\in\L$ and $\omega\in\Omega$:
\be{koz-f}
f_i(\omega_i\mid\omega_{-\infty}^{i-1}) \;=\;
\frac{\displaystyle \sum_{\sigma_{i+1}^0\in\Omega_{i+1}^0}
\exp\biggl[-\sum_{A\in\cS: \atop A\cap[i,0]\neq\emptyset}
\phi_A\big(\sigma_{i+1}^{0}\,\omega_{-\infty}^{i}\big)\biggr]}
{\displaystyle \sum_{\sigma_{i}^0\in\Omega_{i}^0}
\exp\biggl[-\sum_{A\in\cS: \atop A\cap[i,0]\neq\emptyset}
\phi_A\big(\sigma_{i}^{0}\,\omega_{-\infty}^{i-1}\big)\biggr]}
\ee
($\sigma_1^0\equiv \emptyset$).
\ec


\subsection{Application: Long-range Ising ferromagnet chains}
\label{S2.2}

For the alphabet $\cA=\{-1,1\}$, consider the \emph{long range Ising interaction potential} defined by
\be{pot}
\phi_A(\omega)=
\begin{cases}
-\beta\,J(i,j)\,\omega_i\,\omega_j &\text{if } A=\{i,j\},\, i\neq j,\\
0 &\text{otherwise,}
\end{cases}
\ee
with $\sum_j\ |J(i,j)| < \infty$ for each
$i\in\L$ [c.f. \reff{abssum}].  The constants $J(i,j)$ are the \emph{couplings} and $\beta$ is an
overall factor interpreted as the \embf{inverse temperature} (high $\beta$ = low-temperature). The
potential is \embf{ferromagnetic} if $J(i,j)\ge 0$.  Such a potential defines a Gibbsian
\emph{Ising specification} $\gamma^\phi$ through the prescription \reff{gibbsspe}--\reff{gibbsspe*}
and a regular \emph{Ising LSS} $f^\phi$ through \reff{koz-f}.

\bt{ising-res}
Consider an Ising chain and let $J(r):=\sup \{ |J(i,j)|: |i-j|=r\}$.
\begin{itemize}
\item[\rm (a)]
If the chain is ferromagnetic, with decreasing couplings such that
\be{phtr-cond}
\sum_{r\ge 1} \frac{\log\log(r+4)}{r^3 J(r)} \; <\; \infty
\ee
then there are multiple consistent measures at low temperatures and only one at high temperatures
\item[\rm (b)]
If
\be{uniq-cond}
\sum_{j\ge 1} \exp\Bigl[-C\sum_{r\ge 1} (j\wedge r)\, J(r)\Bigr] \;=\;\infty
\quad \forall C>0
\ee
then there is a unique consistent chain at all temperatures.
\end{itemize}
\et
\smallskip

\br{ruel-rem}
Kac and Thompson conjectured (in 1968) that a necessary and sufficient condition for absence of phase
transitions is
\be{kac-tho.for}
\sum_{r\ge 1} r\,J(r) \;<\;\infty
\ee
Part (b) of the precedent theorem shows that the conjecture is false in the half line.  Consider for
instance $J(r)\propto (r^2\log r)^{-1}$.
\er

As an application, consider the  \embf{power-law Ising LSS}:
\be{pli}
J(i,j)\;=\; \frac{1}{|i-j|^p} \quad,\quad p>1\;,
\ee
\bp{th3}
The power law Ising ferromagnet LSS $f$ is well defined if and only $p>1$.  Furthermore:
\begin{itemize}
\item[\rm (i)]
$|\cG(f)|=1$ at high temperature, or at all temperatures if $p> 2$.
\item[\rm (ii)]
If $1<p<2$, then $|\cG(f)|>1$ at low temperature.
\end{itemize}
\ep
The marginal case $p=2$ leads to a very special phase transition in the full line~\cite{frospe82},
and lies outside the scope of our analysis for the half line.

A second application, included for historical reasons, is the \embf{hierarchical Ising chain}.
Its specification version was introduced by Dyson~\cite{dys69a} as a tool for the study
of one-dimensional ferromagnetic Ising models. The hierarchical character of the model
allows for a number of explicit computations that
make its study easier. Moreover, phase transitions for the hierarchical and
the power-law Ising models are related.  The hierarchical chain is defined by considering blocks of
sizes $2^p$, $p\ge 1$ placed consecutively to the left of the origin.  Spins within the same $2^p$-block
interact through a coupling $2^{-2p+1}\,b_p$, for a suitable sequence of positive numbers $b_p$.  Thus,
\be{hier-int}
J(i,j)\;=\; \sum_{n\ge p(i,j)} \frac{b_q}{2^{2q-1}}\;,
\ee
where $p(i,j)$ is the smallest $p$ such that $i$ and $j$ belongs to the same $2^p$ block.

Dyson's results can be transcribed in the following form.
Denote
\be{bpcond}
\Sigma(b)=\sum_{p\geq 1} 2^{-2p+1}\, b_p(2^p-1)
\quad\text{and}\quad
\Sigma^\star(b)=\sum_{p\geq 1}\big(\log(1+p)\big)b_p^{-1}\;.
\ee
\bp{th2}
Assume that $\Sigma(b)<\infty$. Then, the LSS $f$ defined by \reff{koz-f}--\reff{pot} satisfy
\begin{itemize}
\item[{\rm (i)}] $|\cG(f)|=1$ at all temperatures satisfying $\beta\Sigma(b)<1$ or
at all temperatures when $b_p$ are bounded.
\item[{\rm (ii)}] $|\cG(f)|>1$ if $\beta>8 \Sigma^\star(b)$.
\item[{\rm (iii)}] For $b_p=2^{(2-\alpha)p}$ the LSS admits admits several
hierarchical chains at low temperature if and only if $1<\alpha<2$.
\end{itemize}
\ep


\subsection{The uniqueness issue}
\label{S2.3}
Published uniqueness criteria for chains refer to the following four different ways of measuring
the sensitivity of a LSS to changes in the past:
\begin{eqnarray}
\var_k(f_i)
&=&\sup_{\omega\,,\sigma}\Bigl | f_i\bigl(\omega_i\mid\omega_{\le i-1}\bigr)-
f_i\big(\omega_i\mid\omega_{k+1}^{i-1}\,\sigma_{\le k}\big)\Bigr |
\;;\label{var}\\
\osc_k(f_i)
&=&\sup_{\omega\,,\sigma}\Bigl | f_i\bigl(\omega_i\mid\omega_{\le i-1}\bigr)-
f_i\bigl(\omega_i\mid\omega_{k+1}^{i-1}\,\sigma_k\,
\omega_{\le k-1}\bigr)\Bigr |\;;\label{osc}\\
a_k(f_i)&=&\inf_{\sigma} \sum_{\xi_i\in\Omega_i}\inf_{\omega}
f_i\bigl(\xi_i\mid\sigma_k^{i-1}\,\omega_{\le k-1}\,\bigr)\;;\label{const-a}\\
b_k(f_i)&=&\inf_{\omega,\sigma} \sum_{\omega_i\in\cA}
f_i\bigl(\omega_i\mid\omega_{\le i-1}\big)\wedge f_i\big(\omega_i\mid
\omega_{k}^{i-1}\,\sigma_{\le k-1}\bigr) \;.
\label{const-b}
\end{eqnarray}
The first and second quantity are called, respectively, the \embf{$k$-variation} and the
\embf{$k$-oscillation} of the kernel $f_i$.  They are related by the obvious inequalities
\be{varoscine}
\osc_k(f_i) \;\leq\; \var_k(f_i)
\;\leq\; \sum_{j<k}\osc_j(f_i)\;.
\ee
Let us summarize further relations valid for the Ising  case.
\bp{rel-IC}
\begin{itemize}
\item[\rm (a)] If $|\cA|=2$, then for each $k<i\in\L$
\be{a=b}
a_k(f_i)\;=\;b_k(f_i)\;=\;1-\var_k(f_i)\;.
\ee
\item[\rm (b)] If $\{f_i\}$ is an Ising LSS, then for each $k<i\in\L$:
\be{prop-3-0}
\osc_k(f_{i})\leq \beta \sum_{j=i}^0 |J(j,k)|
\quad,\quad
\var_k(f_{i})\leq \beta \sum_{j=i}^0 \sum_{\ell\le k+1}|J(j,\ell)|\;.
\ee
In particular, if $|J(i,j)|\geq |J(k,l)|$ as soon as $|i-j|\leq |k-l|$,
\be{prop-3}
\osc_k(f_{i})\leq \beta \,(i+1)\, |J(i,k)|
\quad,\quad
\var_k(f_{i})\leq \beta \,(i+1) \sum_{j\leq k+1}|J(i,j)|\;.
\ee
\end{itemize}
\ep

The following are the only chain uniqueness criteria proven without a shift-invariance hypothesis.
\bp{unicLSS}
A continuous LSS $f$ admits exactly one consistent chain if it satisfies one of the following assertions:
\begin{itemize}
\item[\rm{(a)}] \emph{CFF \cite{comferfer02}:} $f$ non-null
and
\be{comferfer}
\sum_{j< i}\prod_{k=j}^{i-1} a_k(f_i)=\infty\;,\quad \forall\,i\in\L\;;
\ee
\item[\rm{(b)}] \emph{One-sided boundary-uniformity \cite{fermai05}:}
There exists $C>0$ satisfying:  For every $m\in\Z$ and every cylinder set $A\in\cF_{\leq m}$
there exists $n<m$ such that $\displaystyle f_{[n,m]}(A\mid\xi)\geq C f_{[n,m]}(A\mid \eta)$
for all $\xi,\eta\in\Omega$.  In particular, this condition is satisfied if
\be{osc-nonshift}
\sum_{j< i} \var_j(f_i) \;<\; \infty\;,\quad \forall\,i\in\L\;.
\ee
\end{itemize}
\ep

For comparison purposes, let us list the uniqueness criteria proven for shift-invariant LSS.
\bp{unicgfunct}
A continuous $g$-function $f_0$ admits exactly one consistent chain if it satisfies one of the following criteria:
\begin{itemize}
\item[\rm{(a)}] \emph{Harris \cite{har55,sten03}:} $f_0$ non-null and
$\displaystyle\sum_{j<0}\prod_{k=j}^{-1}\bigg(1-\frac{|\cA|}{2}\var_k(f_0)\bigg)=\infty$.
\item[\rm{(b)}] \emph{Stenflo \cite{sten03}:} $f_0$ non-null and
$\displaystyle\sum_{j< 0}\prod_{k=j}^{-1} b_k(f_0)=\infty$.
\item[\rm{(c)}] \emph{Johansson-\"Oberg \cite{johobe03}:}
$f_0$ non-null and
$\displaystyle\sum_{j< 0}\var_j^2(f_0)<\infty$;
\item[\rm{(d)}] \emph{One-sided Dobrushin \cite{fermai05}:}
$\displaystyle \sum_{j< 0}\osc_j(f_0)\;<\;1$.
\end{itemize}
\ep

It is natural to ask whether these criteria admit a non-shift-invariant version, in which
each $f_0$-condition is replaced by a similar $f_i$-condition valid for all $i\in\L$,
\emph{without asking uniformity with respect to $i$} (that is, without imposing conditions
on $\sup_i f_i$).  For our Ising  examples, (a) of Proposition~\ref{rel-IC} shows that the
non-shift invariant versions of the first two criteria in Proposition~\ref{unicgfunct} coincide
with the CFF criteria.   We are unable to test the corresponding version for one-sided
Dobrushin due to the factor $i+1$ in the leftmost bound in~\reff{prop-3} (though the ``$\sup f_i$" version
of the specification Dobrushin criterium can be applied and yields uniqueness at high
temperature).  On the other hand, the rightmost bound in~\reff{prop-3} allows us to conclude about
the remaining criterion.
\br{johobe-rk}
The power law Ising ferromagnetic LSS $f$ with $3/2<p<2$ show that
the non-shift-invariant version of the Johansson-\"Oberg criterion is false
in general.  Indeed, these LSS do exhibit phase transitions, by Theorem \ref{th3},
but
$\sum_{j\leq i} \var_k^2(f_i)<\infty$ for all $ i\in\Z$.
\er

\section{Proofs}
\label{S3}


\subsection{Proof of Theorem \ref{th1} and Corollary \ref{koz-chains}}
\label{S3.1}

Corollary \ref{koz-chains} is a direct consequence of Theorem \ref{th1} and Kozlov theorem \cite{koz74}.
The proof Theorem \ref{th1} runs as follows.
\smallskip

1)(a) This is a direct consequence of (\ref{1to1-maps**}).
\smallskip

1)(b) Observe that the Gibbsianness of $\gamma^f$ follows directly from the
definition (\ref{1to1-maps*}). Moreover, $\gamma_{\Lambda}^{f}(A\mid \cdot\,)$
is clearly $\mathcal{F}_{\Lambda^{c}}$-measurable for every $\Lambda\in\mathcal{S}$
and every $A \in \mathcal{F}$. Condition (b) of Definition \ref{lis1} together
with the presence of the indicator function $\one\{\omega_{[l_\Lambda,0]\cap\Lambda^\c}\}$
in the numerator of \reff{1to1-maps*} ensure that $\gamma_{\Lambda}^{f}(B \mid \cdot \, )
= \one_B(\, \cdot\,)$ for every $\Lambda \in \mathcal{S}$ and every $B \in
\mathcal{F}_{\Lambda^{c}}$. To conclude the proof that $\gamma^f$ is a Gibbsian specification,
it suffices to show that
\be{th1*-3}
\sum_{\omega_{\Delta \setminus
\Lambda}}\gamma_{\Lambda}^{f}\left( \omega_{\Lambda} \mid
\omega_{\Lambda^{c}} \right) \gamma_{\Delta}^{f}\left( \omega_{\Delta
\setminus \Lambda} \mid \omega_{\Delta^{c}} \right) \;=\;
\gamma_{\Delta}^{f}\left( \omega_{\Lambda} \mid
\omega_{\Delta^{c}}\right)
\ee
for each $\Lambda, \Delta \in \mathcal{S}$ such that $\Lambda \subset
\Delta$ and each $\omega \in \Omega$. Define
$G_\Lambda\colon\cF\ra[0,1]$, $\Lambda\in\cS$, by
\be{Gfunction}
G_\Lambda\big(\,\cdot\mid\omega_{\Lambda^{\c}}\big)
=f_{[l_\Lambda,0]}\Big(\,\cdot\,\one_{\omega_{[l_\Lambda,0]\cap\Lambda^{\c}}}
\mid\omega_{\Lambda_{-}}\Big).
\ee
By \reff{1to1-maps*}
\be{th1*-5}
\gamma_{\Lambda}^f\big(\omega_\Lambda \mid \omega_{\Lambda^{c}}\big)
\;=\;\frac{G_{\Lambda}\left( \one_{\omega_{\Lambda}} \mid
\omega_{\Lambda^{c}}\right)}{G_{\Lambda}\left( \Omega_\Lambda \mid
\omega_{\Lambda^{c}}\right)}\;.
\ee
Using that $f_{[l,n]}=f_{[l,m]}f_{[m+1,n]}$ for all $l\leq m<n\leq 0$,
we obtain
\be{th1*-7}
G_{\Delta}\big(\one_{\omega_{\Delta \setminus \Lambda}}\mid\omega_{\Delta^{c}}\big)
\;f\;_{[l_\Delta,l_{\Lambda}-1]}\Big(\one\big\{\omega_{l_\Delta}^{l_\Lambda-1}\big\}
\mid\omega_{\Delta_{-}}\Big)
\times G_{\Lambda}\big(\Omega_\Lambda\mid\omega_{\Lambda^{c}}\big)
\ee
and
\be{th1*-9}
G_{\Delta}\big(\one_{\omega_{\Delta}}\mid\omega_{\Delta^{c}}\big)\;=\;
f_{[l_\Delta,l_{\Lambda}-1]}\Big(\one\big\{\omega_{l_\Delta}^{l_\Lambda-1}\big\}
\mid\omega_{\Delta_{-}}\Big) \times
G_{\Lambda}\big(\one_{\omega_{\Lambda}}\mid\omega_{\Lambda^{c}}\big)
\;.
\ee
Therefore
\be{th1*-11}
\frac{G_{\Lambda}\big(\one_{\omega_{\Lambda}}\mid \omega_{\Lambda^{c}}\big)}
{G_{\Lambda}\big(\Omega_\Lambda \mid \omega_{\Lambda^{c}}\big)}
\times \frac{G_{\Delta}\big(\one_{\omega_{\Delta \setminus
\Lambda}} \mid \omega_{\Delta^{c}}\big)}{G_{\Delta}\big(\Omega_\Delta \mid
\omega_{\Delta^{c}}\big)}
=\frac{G_{\Delta}\big(\one_{\omega_{\Delta}}\mid\omega_{\Delta^{c}}\big)}
{G_{\Delta}\big(\Omega_\Delta\mid \omega_{\Delta^{c}}\big)}.
\ee
Identity \reff{th1*-3} follows from \reff{th1*-5} and \reff{th1*-11}.
\smallskip

2)(a) For all $\Lambda\in\cS$ and $\omega\in\Omega$
\be{th1-5}
\gamma_\Lambda^{f^\gamma}\big(\omega_\Lambda\mid\omega_{\Lambda^\c}\big)
\;=\;\frac{\displaystyle\gamma_{[l_\Lambda,0]}\big(\omega_{l_\Lambda}^0
\mid\omega_{\Lambda_{-}}\big)}
{\displaystyle\sum_{\sigma_\Lambda\in\Omega_\Lambda}
\gamma_{[l_\Lambda,0]}\big(\sigma_\Lambda\omega_{[l_\Lambda,0]\cap\Lambda^\c}
\mid\omega_{\Lambda_{-}}\big)}\;.
\ee
By the consistency of $\gamma$ [Definition \ref{lis1} (c)]
\be{th1-7}
\gamma_{[l_\Lambda,0]}\big(\omega_{l_\Lambda}^0\mid\omega_{\Lambda_{-}}\big)
\;=\;\gamma_\Lambda\big(\omega_\Lambda\mid\omega_{\Lambda^\c}\big)
\sum_{\sigma_\Lambda\in\Omega_\Lambda}
\gamma_{[l_\Lambda,0]}\big(\sigma_\Lambda\omega_{[l_\Lambda,0]\cap\Lambda^\c}
\mid\omega_{\Lambda_{-}}\big)\;.
\ee
From \reff{th1-5} and \reff{th1-7} we obtain that
$\gamma_\Lambda^{f^\gamma}\big(\omega_\Lambda\mid\omega_{\Lambda^\c}\big)
=\gamma_\Lambda\big(\omega_\Lambda\mid\omega_{\Lambda^\c}\big)$.
\smallskip

2)(b) For all $\Lambda\in\cS_b$ and $\omega\in\Omega$
\be{th1-9}
f_\Lambda^{\gamma^f}\big(\omega_\Lambda\mid\omega_{\Lambda_{-}}\big)
\;=\;\gamma_{[l_\Lambda,0]}^f\big(\omega_\Lambda\mid\omega_{\Lambda_{-}}\big)
\;=\;\frac{f_{[l_\Lambda,0]}\big(\omega_\Lambda\mid\omega_{\Lambda_{-}}\big)}
{f_{[l_\Lambda,0]}\big(\Omega_{l_\Lambda}^0\mid\omega_{\Lambda_{-}}\big)}
\;=\;f_\Lambda\big(\omega_\Lambda\mid\omega_{\Lambda_{-}}\big).\;
\ee
\smallskip

2)(c) We use \reff{2.4.1} and consistency to obtain the following two strings of inequalities.
If $\mu\in\cG(\gamma)$ and $\Lambda\in\cS_b$,
\be{th1-11}
\mu \,f_\Lambda^\gamma
\;=\;\mu \,\gamma_{[l_\Lambda,0]}
\;=\;\mu\;.
\ee
If $\mu\in\cG\big(f^\gamma\big)$ and $\Lambda\in\cS$,
\be{th1-13}
\mu\,\gamma_\Lambda
=\bigl(\mu \,f_{[l_\Lambda,0]}^{\gamma}\bigr)\gamma_\Lambda
=\mu \bigl(\gamma_{[l_\Lambda,0]}\,\gamma_\Lambda\bigr)
=\mu \,\gamma_{[l_\Lambda,0]}
=\mu \,f_{[l_\Lambda,0]}^{\gamma}
=\mu\;.
\ee
Together, \reff{th1-11} and \reff{th1-13} show that $\cG\big(f^\gamma\big)=\cG(\gamma)$.
\smallskip

2)(d) The following two displays are a consequence of \reff{2.2.1} and consistency.
For any $\mu\in\cG(f)$ and $\Lambda\in\cS$,
\be{th1-15}
\mu\,\gamma_\Lambda^f
\;=\;\bigl(\mu\, f_{[l_\Lambda,0]}\bigr)\gamma_\Lambda^f
\;=\;\mu\bigl (\gamma_{[l_\Lambda,0]}^f\,\gamma_\Lambda^f\bigr)
\;=\; \mu\,\gamma_{[l_\Lambda,0]}^f
\;=\; \mu\,f_{[l_\Lambda,0]}
\;=\;\mu\;.
\ee
For any $\mu\in\cG\big(\gamma^f\big)$ and $\Lambda\in\cS_b$,
\be{th1-17}
\mu\, f_\Lambda
\;=\;\Big(\mu\,\gamma_{[l_\Lambda,0]}^f\Big)f_\Lambda
\;=\;\mu\bigl(f_{[l_\Lambda,0]}\,f_\Lambda\bigr)
\;=\;\mu\, f_{[l_\Lambda,0]}
\;=\;\mu\,\gamma_{[l_\Lambda,0]}^f
\;=\;\mu\;.
\ee
The combination of both lines shows that $\cG\big(\gamma^f\big)=\cG(f)$.
\qed

\subsection{Proof of Proposition \ref{rel-IC}}
\label{S3.1.1}

We need the following well known bound, whose proof we present for completeness
(we follow the approach of~\cite[Lemma V.1.4]{sim93}).
\bl{simlem}
Let $\gamma^\varphi$ be a Gibbsian specification for some absolutely summable
potential $\varphi$. Then, for any $\Lambda\in\cS$, $h\in\cF_\Lambda$ such that
$\|h\|_\infty<\infty$ and $\omega,\sigma\in\Omega$,
\be{simlem-3}
\bigg|\int_{\Omega_\Lambda}h(\xi)\Big(
\gamma_\Lambda^\varphi\big(d\xi\mid\omega_{\Lambda^\c}\big)
-\gamma_\Lambda^\varphi\big(d\xi\mid\sigma_{\Lambda^\c}\big)\Big)\bigg|
\leq \|h\|_\infty \sup_{\xi\in\Omega_\Lambda}
\Big|H_{\Lambda,\omega}^\varphi(\xi)
-H_{\Lambda,\sigma}^\varphi(\xi)\Big|,
\ee
where $H^\varphi$ is the Hamiltonian associated to $\varphi$.
\el

\bpr
For all $\Lambda\in\cS$, $\omega,\sigma,\xi\in\Omega$ and $0<\theta<1$, define
$\Gamma_{\Lambda,\omega,\sigma}^{\varphi,\theta}\colon\Omega_\Lambda\ra(0,1)$ by
\be{simlem-11}
\Gamma_{\Lambda,\omega,\sigma}^{\varphi,\theta}(\xi)
=\frac{\displaystyle\exp\Big[\theta H_{\Lambda,\omega}^{\varphi}(\xi)
+(1-\theta)H_{\Lambda,\sigma}^{\varphi}(\xi)\Big]}
{\displaystyle\sum_{\eta\in\Omega_\Lambda}
\exp\Big[\theta H_{\Lambda,\omega}^{\varphi}(\eta)
+(1-\theta)H_{\Lambda,\sigma}^{\varphi}(\eta)\Big]}.
\ee
Then, to prove (\ref{simlem-3}), it suffices to see that
\be{simlem-13}
\begin{aligned}
&\bigg|\int_{\Omega_\Lambda}h(\xi)\Big(
\gamma_\Lambda^\varphi\big(d\xi\mid\omega_{\Lambda^\c}\big)
-\gamma_\Lambda^\varphi\big(d\xi\mid\sigma_{\Lambda^\c}\big)\Big)\bigg|\\
&\quad\leq \int_{0}^{1}\bigg|\frac{d}{d\theta}\bigg[\int_{\Omega_\Lambda} h\,
d\Gamma_{\Lambda,\omega,\sigma}^{\varphi,\theta}\bigg]\bigg|\, d\theta\\
&\quad=\int_{0}^{1}\bigg|\int_{\Omega_\Lambda} h\Big(H_{\Lambda,\omega}^{\varphi}
-H_{\Lambda,\sigma}^{\varphi}\Big) d\Gamma_{\Lambda,\omega,\sigma}^{\varphi,\theta}
-\int_{\Omega_\Lambda} h\, d\Gamma_{\Lambda,\omega,\sigma}^{\varphi,\theta}
\int_{\Omega_\Lambda} \Big(H_{\Lambda,\omega}^{\varphi}
-H_{\Lambda,\sigma}^{\varphi}\Big) d\Gamma_{\Lambda,\omega,\sigma}^{\varphi,\theta}
\bigg|\, d\theta\\
&\quad\leq \|h\|_\infty\int_{0}^{1}\int_{\Omega_\Lambda} \bigg|\Big(H_{\Lambda,\omega}^{\varphi}
-H_{\Lambda,\sigma}^{\varphi}\Big) -\int_{\Omega_\Lambda} \Big(H_{\Lambda,\omega}^{\varphi}
-H_{\Lambda,\sigma}^{\varphi}\Big)d\Gamma_{\Lambda,\omega,\sigma}^{\varphi,\theta} \bigg|\,
d\Gamma_{\Lambda,\omega,\sigma}^{\varphi,\theta}\, d\theta\\
&\quad\leq \|h\|_\infty\sup_{\xi_\Lambda\in\Omega_\Lambda}\big|H_{\Lambda,\omega}^{\varphi}
-H_{\Lambda,\sigma}^{\varphi}\big|.
\end{aligned}
\ee
\epr

\noindent
\emph{Proof of Proposition \ref{rel-IC}}.
\smallskip

(a) For any $k<i$,
\be{th3-60}
\begin{aligned}
a_k(f_i)
&=\inf_{\sigma,\omega,\xi}\Bigl[
f_i(1\mid\sigma_k^{i-1}\,\omega_{\le k-1})
+f_i(-1\mid\sigma_k^{i-1}\,\xi_{\le k-1})\Bigr]\\
&=1-\sup_{\sigma,\omega,\xi}\Bigl[
-f_i(1\mid\sigma_k^{i-1}\,\omega_{\le k-1})
+f_i(1\mid\sigma_k^{i-1}\,\xi_{\le k-1})\Bigr]\\
&=1-\var_k(f_i).
\end{aligned}
\ee
Likewise,
\begin{eqnarray}
\label{th3-60.1}
b_k(f_i)&=&\inf_{\sigma,\omega,\xi}\Bigl[
f_i(1\mid\sigma_k^{i-1}\,\omega_{\le k-1})\wedge
f_i(1\mid\sigma_k^{i-1}\,\xi_{\le k-1})\nonumber\\
&&\qquad {}+\bigl[1-f_i(1\mid\sigma_k^{i-1}\,\omega_{\le k-1})\bigr] \wedge
\bigl[1 - f_i(1\mid\sigma_k^{i-1}\,\xi_{\le k-1})\bigr]
\Bigr]\nonumber\\
&=&1-\sup_{\sigma,\omega,\xi}\Bigl[
f_i(1\mid\sigma_k^{i-1}\,\omega_{\le k-1})\vee
f_i(1\mid\sigma_k^{i-1}\,\xi_{\le k-1})\\
&&\qquad\qquad {}- f_i(1\mid\sigma_k^{i-1}\,\omega_{\le k-1})\wedge
f_i(1\mid\sigma_k^{i-1}\,\xi_{\le k-1})\Bigr]\nonumber\\
&=&1-\var_k(f_i)\;.\nonumber
\end{eqnarray}
\smallskip

(b) Applying the previous lemma for $\Lambda=[i,0]$ and $h(\xi)=1$ if $\xi_i=\omega_i$
and 0 otherwise, we obtain that, for any $\omega_{\le i-1}, \sigma_{\le i-1}\in \Omega_{\le i-1}$,
\be{simlem-5}
\begin{aligned}
\Bigl|f_i\bigl(\omega_i\mid\omega_{\le i-1}\bigr)
-f_i\bigl(\omega_i\mid\sigma_{\le i-1}\bigr)\Bigr|
\;&=\; \Bigl|\gamma^\phi_{[i,0]}(h\mid\omega_{\le i-1}) -
\gamma^\phi_{[i,0]}(h\mid\sigma_{\le i-1})\Bigr|\\
\;&\le\;
\sup_{\xi}\Bigl|H_{[i,0]}^\phi\bigl(\xi_{i}^{0}\mid\omega_{\le i-1}\big)
-H_{[i,0]}^\phi\big(\xi_{i}^{0}\mid\sigma_{\le i-1}\bigr)\Bigr|\;.
\end{aligned}
\ee
Both inequalities in \reff{prop-3-0} are an immediate consequence.  \qed

\subsection{Proof of Theorem \ref{ising-res}}
\label{S3.2}

(a) This is a direct consequence of Theorem \ref{th1} and Theorem 1 in \cite{dys69a}.
\smallskip

(b) We show the validity of the CFF condition [Proposition \ref{unicLSS}(a)].
Pick an  $\alpha>1$.  If $x$ is small enough, $1-x\ge \exp(-\alpha\,x)$.  Hence,
there exists $j_0<i$ small enough and $K>0$ such that
\be{th3-61}
\sum_{j<i}\prod_{k=j}^{i-1}a_k(f_i)
\;\geq\; K\sum_{j\leq j_0}\exp\Bigg[-\alpha\sum_{k=j}^{j_0}\var_k(f_i)\Bigg].
\ee
Then, by (\ref{uniq-cond}) and (\ref{prop-3}),
\begin{eqnarray}
\label{th3-63}
\sum_{j<i}\prod_{k=j}^{i-1}a_k(f_i)
&\geq&K\sum_{j\leq j_0}\exp\Bigg[-\alpha\beta K(i)\sum_{k=j}^{j_0}\, \sum_{r\geq |k|}J(r)\Bigg]
\nonumber\\
&\geq& K\sum_{j\geq|j_0|}\exp\Bigg[-\alpha\beta K(i)\sum_{r\geq |j_0|}(j\wedge r)J(r)\Bigg]\\
&=&\infty
\end{eqnarray}
by hypothesis. \qed


\end{document}